\documentclass[11pt,a4paper]{article}
\usepackage[utf8]{inputenc}
\usepackage[T1]{fontenc}
\usepackage[english]{babel}
\usepackage{geometry}
\usepackage{amsmath,amssymb,amsthm,mathtools}
\usepackage{newtxtext}
\usepackage[subscriptcorrection]{newtxmath}

\usepackage{booktabs}
\usepackage{cite}
\usepackage[colorlinks,citecolor=blue,urlcolor=blue,linkcolor=red]{hyperref}
\usepackage[capitalize,noabbrev]{cleveref}
\usepackage{doi}

\usepackage{verbatim}
\usepackage{comment}

\numberwithin{equation}{section}
\theoremstyle{plain}
\newtheorem{theorem}{Theorem}[section]
\newtheorem{proposition}[theorem]{Proposition}
\newtheorem{lemma}[theorem]{Lemma}
\newtheorem{corollary}[theorem]{Corollary}
\theoremstyle{remark}

\title{The multivariate Berry–Esseen theorem with an absolute constant}
\author{Zhi Ling\thanks{National University of Singapore. Email: \texttt{lingzhi@nus.edu.sg}.}}
\date{}

\begin{document}
\maketitle

\begin{abstract}
We prove that the multivariate Berry--Esseen inequality for sums of independent random vectors holds with the same constant in every dimension. For independent centered random vectors $X_1,\ldots,X_n\in\mathbb{R}^d$ whose covariance matrices sum to $I_d$, we show that the distance between the distribution of their sum and $N(0,I_d)$, measured uniformly over convex sets, is at most $C_0\sum_{i=1}^n\mathbb{E}\|X_i\|_2^3$ for a universal constant $C_0$. The bound requires only finite third moments and holds uniformly in the dimension, the number of summands, and their distributions. It also yields the corresponding bound for independent and identically distributed vectors. In particular, $C_0$ can be chosen to be $980$.
\end{abstract}

\section{Introduction}

The multivariate Berry--Esseen problem asks how accurately a sum of independent random vectors can be approximated by a Gaussian vector. We investigate bounds on the Gaussian approximation error over convex sets in terms of the sum of the third absolute moments of the summands. The question is whether the multiplicative constant can be chosen independently of the dimension.

Fix integers $n,d\ge1$, and let $\mathcal C_d:=\{C\subseteq\mathbb{R}^d:C\text{ is convex and Borel}\}$. Let $X_1,\ldots,X_n$ be independent random vectors such that
\begin{equation}
\label{eq.normalization}
\mathbb{E}X_i=0,\qquad
\sum_{i=1}^n\operatorname{Cov}(X_i)=I_d,\qquad
\sum_{i=1}^n\mathbb{E}\|X_i\|_2^3<\infty.
\end{equation}
Set $S:=\sum_{i=1}^nX_i$, let $Z\sim N(0,I_d)$, and define
\begin{equation}
\Delta_{n,d}:=\sup_{C\in\mathcal C_d}
\bigl|\mathbb{P}(S\in C)-\mathbb{P}(Z\in C)\bigr|,
\qquad
\beta_{n,d}:=\sum_{i=1}^n\mathbb{E}\|X_i\|_2^3.
\end{equation}
The optimal constant for independent summands is
\begin{equation}
K_d:=\sup\left\{\frac{\Delta_{n,d}}{\beta_{n,d}}:
n\ge1,\ X_1,\ldots,X_n\text{ satisfy \eqref{eq.normalization}}\right\}.
\end{equation}
The covariance normalization ensures that $\beta_{n,d}>0$.

\begin{theorem}[Berry--Esseen bound for independent sums]
\label{thm.main_result}
There is a universal constant $C_0<\infty$ such that $K_d\le C_0$ for every $d\ge1$. Equivalently, any independent random vectors $X_1,\ldots,X_n$ satisfying \eqref{eq.normalization} obey
$$
\Delta_{n,d}\le C_0\beta_{n,d}.
$$
In particular, one may take $C_0=980$.
\end{theorem}

\begin{corollary}[Berry--Esseen bound for i.i.d. sums]
Let $X_1,\ldots,X_n$ be independent and identically distributed random vectors in $\mathbb{R}^d$ satisfying
$$
\mathbb{E}X_1=0,\qquad
\operatorname{Cov}(X_1)=I_d,\qquad
\mathbb{E}\|X_1\|_2^3<\infty.
$$
Let $Z\sim N(0,I_d)$, and let $\mathcal C_d$ be the collection of all convex Borel subsets of $\mathbb{R}^d$. There is a universal constant $C_0<\infty$ such that
$$
\sup_{C\in\mathcal C_d}
\left|\mathbb{P}\left(\frac1{\sqrt n}\sum_{i=1}^nX_i\in C\right)
-\mathbb{P}(Z\in C)\right|
\le\frac{C_0\,\mathbb{E}\|X_1\|_2^3}{\sqrt n}.
$$
\end{corollary}

\subsection{Background and related work}

The classical Berry--Esseen theorem gives a quantitative bound on the error in the central limit theorem in terms of the third absolute moment. Berry \cite{Berry1941} and Esseen \cite{Esseen1945} established the univariate theory using characteristic functions. In the multivariate case, a central question is whether the coefficient $K_d$ in the third-moment bound depends on the dimension $d$.

For sums of independent and identically distributed random vectors, Nagaev \cite{Nagaev1976} obtained a coefficient of order $d$ using Fourier methods. Further bounds with linear dependence on the dimension were obtained by Senatov \cite{Senatov1980} and Sazonov \cite{Sazonov1981}. Bentkus \cite{Bentkus1986} improved the dependence to $d^{1/2}$ using the Lindeberg--Bergstr\"om method. He subsequently obtained the explicit coefficient $400d^{1/4}$ in the i.i.d. case \cite{Bentkus2003} and extended the order $d^{1/4}$ to independent, not necessarily identically distributed summands \cite{Bentkus2005}.

Stein's method provides an alternative approach to this problem \cite{Barbour1990}. G\"otze \cite{Gotze1991} combined the multivariate Stein equation with induction to obtain a coefficient of order $d$. Chen and Fang \cite{ChenFang2015} developed a concentration inequality approach and obtained the explicit coefficient $115d^{1/2}$ for independent sums. Their method controls the probability that a sum lies near the boundary of a convex set. Combining Stein's method with arguments of Bentkus, Rai\v c \cite[Theorem~1.1]{Raic2019} obtained
$$
K_d\le42d^{1/4}+16.
$$

The dimension dependence in these bounds is closely related to Gaussian surface area. Smoothing the indicator function of a convex set introduces an error controlled by the Gaussian measure of a neighborhood of its boundary. Ball \cite{Ball1993} proved that the Gaussian surface area of a convex set is at most $4d^{1/4}$, and Nazarov \cite{Nazarov2003} showed that the order $d^{1/4}$ is optimal. However, the optimal order of Gaussian surface area does not determine the optimal order of the Berry--Esseen coefficient, as noted by Rai\v c \cite{Raic2019}. Bentkus \cite{Bentkus2003} explicitly raised the question of whether the dimension factor could be removed.

Theorem~\ref{thm.main_result} removes this factor under finite third moments. Together with Nagaev's lower-bound construction \cite{Nagaev1976}, our result shows that $K_d$ is of constant order. Our proof combines Stein's method with smoothing scales adapted to the covariance structure in different directions. We obtain the required derivative estimates by replacing individual summands with independent copies. The bound holds uniformly over all dimensions and all independent sums satisfying \eqref{eq.normalization}.

\subsection{Organization of the paper}

Section~\ref{sec.covariance_layers} develops the covariance estimates used to choose different smoothing scales in different directions. Section~\ref{sec.smoothing} constructs smooth approximations to indicator functions of convex sets and controls the resulting approximation error. Sections~\ref{sec.resampling} and~\ref{sec.stein_integral} establish the derivative bounds needed in Stein's method by replacing individual summands with independent copies. Section~\ref{sec.main_proof} combines these estimates to prove the main theorem by induction on $n$. Appendix~\ref{app.explicit_constant} verifies the explicit choice of $C_0$.

\subsection{Notation}

For a vector, $\|\cdot\|_2$ denotes the Euclidean norm. For a matrix, $\|\cdot\|_{\mathrm{op}}$ and $\|\cdot\|_{\mathrm{HS}}$ denote the operator and Hilbert--Schmidt norms. We write $A\succeq B$ when $A-B$ is positive semidefinite, and $A^\dagger$ denotes the Moore--Penrose inverse of a positive semidefinite matrix. The notation $D^mf(x)[v_1,\ldots,v_m]$ denotes a directional derivative, and $D_vf:=Df[v]$. Let $\gamma_d$ denote standard Gaussian measure and $B_2^d:=\{x\in\mathbb{R}^d:\|x\|_2\le1\}$. For random vectors $U,V$ in $\mathbb{R}^d$, define
$$
d_{\mathcal C}(U,V):=\sup_{C\in\mathcal C_d}
|\mathbb{P}(U\in C)-\mathbb{P}(V\in C)|.
$$
We use the same notation when one argument above is a probability distribution. For a fixed array, abbreviate $\beta:=\beta_{n,d}$. The letter $C$ denotes a convex set, and $c$ denotes a positive universal constant that may change from line to line. Unless otherwise stated, all copies and auxiliary Gaussian vectors are independent of $X_1,\ldots,X_n$.

\section{Covariance blocks and nested subspaces}
\label{sec.covariance_layers}

In this section, we partition the summands into three independent blocks and construct nested subspaces to guide the choice of smoothing widths in Section~\ref{sec.smoothing}.

Set $\Sigma_i:=\operatorname{Cov}(X_i)$ and $b_i:=\mathbb{E}\|X_i\|_2^3$. Fix universal parameters $b,\lambda,r_0,q$ such that
$$
0<b<\frac13,\qquad 0<\lambda<\frac13,\qquad
r_0>0,\qquad \lambda+\frac4{r_0}<\frac13,\qquad q>1.
$$
The parameters $b$ and $\lambda$ are covariance thresholds, $r_0$ sets the largest resampling radius, and $q$ is the ratio between successive radii. Choose $\beta_0\in(0,1)$ sufficiently small that
\begin{equation}
\label{eq.small_budget_conditions}
\sqrt{\frac23}\,\beta_0^{2/3}\le\frac13-\lambda-\frac4{r_0},
\qquad
\left(1+\sqrt{\frac23}\right)\beta_0^{2/3}\le\frac13-b.
\end{equation}
These parameters remain fixed throughout the proof. Universal constants may depend on this fixed choice, but not on $n,d,\beta$, or the smoothing parameter introduced below. Throughout this section assume $\beta\le\beta_0$.

\begin{lemma}[Covariance Blocks]
There is a deterministic partition $I_1,I_2,I_3$ of $\{1,\ldots,n\}$ such that, with $B_a:=\sum_{i\in I_a}\Sigma_i$,
\begin{equation}
\label{eq.block_covariance}
B_a\succeq\left(\lambda+\frac4{r_0}\right)I_d
\qquad(a=1,2,3).
\end{equation}
For every $i\in I_a$, one also has $B_a-\Sigma_i\succeq bI_d$.
\end{lemma}

\begin{proof}
Positivity and Lyapunov's inequality give
$$
\begin{aligned}
\sum_i\|\Sigma_i\|_{\mathrm{HS}}^2
&\le\sum_i(\operatorname{tr}\Sigma_i)^2
=\sum_i(\mathbb{E}\|X_i\|_2^2)^2\\
&\le\sum_i b_i^{4/3}
\le\left(\sum_i b_i\right)^{4/3}
=\beta^{4/3}.
\end{aligned}
$$
Assign each index independently and uniformly to one of three blocks. Since $\sum_i\Sigma_i=I_d$,
$$
\begin{aligned}
\mathbb{E}\sum_{a=1}^3\left\|B_a-\frac13I_d\right\|_{\mathrm{HS}}^2
&=\sum_{a=1}^3\sum_i
\operatorname{Var}(\mathbf1_{\{i\in I_a\}})\|\Sigma_i\|_{\mathrm{HS}}^2\\
&=\frac23\sum_i\|\Sigma_i\|_{\mathrm{HS}}^2
\le\frac23\beta^{4/3}.
\end{aligned}
$$
Choose a partition for which the sum does not exceed its expectation. Then
$$
\left\|B_a-\frac13I_d\right\|_{\mathrm{op}}
\le\sqrt{\frac23}\,\beta^{2/3}.
$$
Also,
$$
\|\Sigma_i\|_{\mathrm{op}}
\le\operatorname{tr}\Sigma_i
=\mathbb{E}\|X_i\|_2^2
\le b_i^{2/3}\le\beta^{2/3}.
$$
The two assertions now follow from \eqref{eq.small_budget_conditions}:
$$
\begin{aligned}
B_a&\succeq\left(\frac13-\sqrt{\frac23}\,\beta^{2/3}\right)I_d
\succeq\left(\lambda+\frac4{r_0}\right)I_d,\\
B_a-\Sigma_i
&\succeq\left[\frac13-\left(1+\sqrt{\frac23}\right)\beta^{2/3}\right]I_d
\succeq bI_d.
\end{aligned}
$$
\end{proof}

Fix this partition. Let $\widehat X_i$ be an independent copy of $X_i$, and put $D_i:=\widehat X_i-X_i$. For $R>0$, define
$$
\Gamma_a(R):=\frac12\sum_{i\in I_a}
\mathbb{E}\bigl[D_iD_i^\top\mathbf1_{\{\|D_i\|_2\le R\}}\bigr],
\qquad
Q_a(R):=B_a-\Gamma_a(R).
$$
Since $\mathbb{E}[D_iD_i^\top]=2\Sigma_i$, the deficit has the exact representation
\begin{equation}
\label{eq.deficit_identity}
Q_a(R)=\frac12\sum_{i\in I_a}
\mathbb{E}\bigl[D_iD_i^\top\mathbf1_{\{\|D_i\|_2>R\}}\bigr]\succeq0.
\end{equation}
Moreover,
\begin{equation}
\label{eq.displacement_moment}
\sum_i\mathbb{E}\|D_i\|_2^3
\le4\sum_i\bigl(\mathbb{E}\|\widehat X_i\|_2^3+\mathbb{E}\|X_i\|_2^3\bigr)
=8\beta.
\end{equation}

\begin{lemma}[Covariance Deficits]
Set $R_j:=r_0\beta q^j$ for integers $j\le0$. Then
\begin{equation}
\label{eq.deficit_budget}
\sum_a\operatorname{tr}Q_a(R)\le\frac{4\beta}{R},
\qquad
\sum_{j\le0}R_j\sum_a\operatorname{tr}Q_a(R_j)
\le\frac{4q}{q-1}\beta.
\end{equation}
In particular, $\Gamma_a(R_0)\succeq\lambda I_d$ for every block.
\end{lemma}

\begin{proof}
Taking traces in \eqref{eq.deficit_identity} gives
$$
\sum_a\operatorname{tr}Q_a(R)
=\frac12\sum_i\mathbb{E}\bigl[\|D_i\|_2^2\mathbf1_{\{\|D_i\|_2>R\}}\bigr]
\le\frac1{2R}\sum_i\mathbb{E}\|D_i\|_2^3
\le\frac{4\beta}{R}.
$$
For every $s>0$, the geometric progression satisfies
$$
\sum_{j\le0:R_j<s}R_j\le\frac{q}{q-1}s.
$$
Consequently, Tonelli's theorem yields
$$
\begin{aligned}
\sum_{j\le0}R_j\sum_a\operatorname{tr}Q_a(R_j)
&=\frac12\sum_i\mathbb{E}\left[\|D_i\|_2^2
\sum_{j\le0:R_j<\|D_i\|_2}R_j\right]\\
&\le\frac{q}{2(q-1)}\sum_i\mathbb{E}\|D_i\|_2^3
\le\frac{4q}{q-1}\beta.
\end{aligned}
$$
Finally,
$$
\Gamma_a(R_0)
=B_a-Q_a(R_0)
\succeq\left(\lambda+\frac4{r_0}-\frac{4\beta}{r_0\beta}\right)I_d
=\lambda I_d.
$$
\end{proof}

Let $G_{a,j}$ denote the spectral subspace of $\Gamma_a(R_j)$ corresponding to eigenvalues at least $\lambda$. Define
\begin{equation}
\label{eq.nested_subspaces}
E_j:=\bigcap_{a=1}^3\bigcap_{\ell=j}^0G_{a,\ell},
\qquad
H_j:=E_j\cap E_{j-1}^{\perp}.
\end{equation}
Although $\Gamma_a(R_j)$ increases with $j$ in the positive semidefinite order, its spectral subspaces $G_{a,j}$ need not be nested. The intersection across scales enforces nesting while preserving the spectral lower bound for every block.

Let $\Pi_j$ be the orthogonal projection onto $E_j$, and set $P_j:=\Pi_j-\Pi_{j-1}$. Thus $P_j$ projects onto $H_j$, $E_j\subseteq E_{j+1}$, and $E_0=\mathbb{R}^d$. If $N_a:=|I_a|$, then
$$
\operatorname{tr}\Gamma_a(R_j)
\le\frac{N_aR_j^2}{2}\longrightarrow0
\qquad(j\to-\infty).
$$
Hence $E_j=\{0\}$ for all sufficiently negative $j$. Only finitely many $H_j$ are nonzero, and $\sum_{j\le0}P_j=I_d$.

\begin{lemma}[Spectral Layers]
For $u\in E_j$ and every block,
\begin{equation}
\label{eq.inverse_energy}
\Gamma_a(R_j)\Gamma_a(R_j)^\dagger u=u,
\qquad
u^\top\Gamma_a(R_j)^\dagger u\le\lambda^{-1}\|u\|_2^2.
\end{equation}
The layers satisfy
\begin{equation}
\label{eq.dimension_budget}
\sum_{j\le0}q^j\dim H_j\le\frac{q^2}{q-1}.
\end{equation}
\end{lemma}

\begin{proof}
The first assertions follow from $E_j\subseteq G_{a,j}$. Let $P_{a,j}^{\mathrm{bad}}$ be the spectral projection onto eigenvalues below $\lambda$, and let $s_{a,j}$ be its rank. By \eqref{eq.block_covariance},
$$
\begin{aligned}
\operatorname{tr}Q_a(R_j)
&\ge\operatorname{tr}\bigl(P_{a,j}^{\mathrm{bad}}Q_a(R_j)\bigr)\\
&=\operatorname{tr}\bigl(P_{a,j}^{\mathrm{bad}}B_a\bigr)
-\operatorname{tr}\bigl(P_{a,j}^{\mathrm{bad}}\Gamma_a(R_j)\bigr)\\
&\ge\left(\lambda+\frac4{r_0}-\lambda\right)s_{a,j}
=\frac4{r_0} s_{a,j}.
\end{aligned}
$$
The codimension of an intersection is at most the sum of its codimensions. Writing $q_j:=\operatorname{codim}E_j$ and using $G_{a,0}=\mathbb{R}^d$, we obtain
$$
q_j\le\frac{r_0}4\sum_{\ell=j}^{-1}\sum_a\operatorname{tr}Q_a(R_\ell).
$$
Since $q_0=0$ and $R_jq_j\to0$ as $j\to-\infty$, summation by parts and \eqref{eq.deficit_budget} give
$$
\begin{aligned}
\sum_{j\le0}R_j\dim H_j
&=\sum_{j\le0}R_j(q_{j-1}-q_j)
=(q-1)\sum_{j<0}R_jq_j\\
&\le\frac{r_0(q-1)}4\sum_{\ell<0}
\left(\sum_{j\le\ell}R_j\right)\sum_a\operatorname{tr}Q_a(R_\ell)\\
&=\frac{r_0q}{4}\sum_{\ell<0}R_\ell\sum_a\operatorname{tr}Q_a(R_\ell)
\le\frac{r_0q^2}{q-1}\beta.
\end{aligned}
$$
Dividing by $r_0\beta$ proves \eqref{eq.dimension_budget}.
\end{proof}

\section{Directional smoothing and Gaussian comparison}
\label{sec.smoothing}

Recall that the orthogonal layers $H_j$ satisfy the weighted dimension bound \eqref{eq.dimension_budget}. We use this bound to construct smooth approximations to convex indicators and then establish the Gaussian comparison estimates needed for the derivative bounds.

For a subspace $H$, define $B_H:=\{x\in H:\|x\|_2\le1\}$. For a closed convex set $C$ and a compact symmetric convex set $F$, write
$$
C+F:=\{x+y:x\in C,\ y\in F\},
\qquad
C\ominus F:=\{x:x+F\subseteq C\}.
$$

\begin{lemma}[Gaussian Sections]
\label{lem.gaussian_sections}
Let $H_1,\ldots,H_m$ be mutually orthogonal nonzero subspaces with sum $\mathbb{R}^d$, and let $P_1,\ldots,P_m$ be their projections. If $M:=\sum_{j=1}^m\delta_jP_j$ with $\delta_j>0$, then every closed convex set $C$ satisfies
\begin{equation}
\label{eq.section_bound}
\max\left\{\gamma_d(C+MB_2^d)-\gamma_d(C),\
\gamma_d(C)-\gamma_d(C\ominus MB_2^d)\right\}
\le c\sum_{j=1}^m\delta_j(\dim H_j)^{1/4}.
\end{equation}
For a Gaussian vector whose covariance is at least $sI_d$, the bound holds uniformly in its mean after multiplication by $s^{-1/2}$.
\end{lemma}

\begin{proof}
The Gaussian perimeter bound and boundary-neighborhood identity in \cite[Theorem~1.2 and Proposition~1.1]{Raic2019} imply, in dimension $r\ge1$,
$$
\max\{\gamma_r(A+tB_2^r)-\gamma_r(A),
\gamma_r(A)-\gamma_r(A\ominus tB_2^r)\}
\le c\,tr^{1/4}.
$$
Fix a subspace $H$ of dimension $r$. A standard Gaussian vector splits into independent standard Gaussian vectors on $H$ and $H^\perp$. For $y\in H^\perp$, put $C_y:=\{x\in H:x+y\in C\}$. Since
$$
(C+tB_H)_y=C_y+tB_H,\qquad
(C\ominus tB_H)_y=C_y\ominus tB_H,
$$
integrating the $r$-dimensional bound over $y$ gives the same estimate for expansion and erosion along $H$.

Set $F:=\sum_j\delta_jB_{H_j}$. The inclusion $MB_2^d\subseteq F$ follows from $\|P_jx\|_2\le1$ for $\|x\|_2\le1$. Successively adding the sets $\delta_jB_{H_j}$ and summing the section estimates bounds the outer difference in \eqref{eq.section_bound}. For the inner difference, use $C\ominus MB_2^d\supseteq C\ominus F$ and
$$
C\ominus F
=\bigl(\cdots((C\ominus\delta_1B_{H_1})\ominus\delta_2B_{H_2})\cdots\bigr)
\ominus\delta_mB_{H_m}.
$$
The same section estimate applies at each step, including when an erosion is empty.

Finally, a Gaussian vector with mean $\mu$ and covariance $V\succeq sI_d$ can be written as $\mu+\sqrt sZ+Z'$, where $Z'$ is independent Gaussian with covariance $V-sI_d$. Conditional on $Z'$, translate $C$ and divide by $\sqrt s$. Apply \eqref{eq.section_bound} and average.
\end{proof}

For the layers in \eqref{eq.nested_subspaces}, fix $L\ge1$ and set
\begin{equation}
\delta_j:=L\beta\,q^{j/3},
\qquad
M:=\sum_{j\le0}\delta_jP_j.
\end{equation}
This is a finite sum over the nonzero layers, so $M$ is positive definite. By H\"older's inequality and \eqref{eq.dimension_budget},
\begin{equation}
\label{eq.smoothing_cost}
\begin{aligned}
\sum_{j\le0}\delta_j(\dim H_j)^{1/4}
&=L\beta\sum_{j\le0}(q^j\dim H_j)^{1/4}q^{j/12}\\
&\le L\beta
\left(\sum_{j\le0}q^j\dim H_j\right)^{1/4}
\left(\sum_{j\le0}q^{j/9}\right)^{3/4}\\
&\le\frac{(q^2/(q-1))^{1/4}}{(1-q^{-1/9})^{3/4}}L\beta.
\end{aligned}
\end{equation}
Combining this estimate with Lemma~\ref{lem.gaussian_sections}, fix a universal constant $c_{\mathrm s}$ such that both Gaussian boundary differences for $MB_2^d$ are bounded by $c_{\mathrm s}L\beta$. This constant is fixed independently of $L$.

For nonempty closed convex $C$, define
\begin{equation}
h_C(x):=\psi\bigl(\operatorname{dist}(M^{-1}x,M^{-1}C)\bigr),
\qquad
\psi(t):=
\begin{cases}
1-2t^2,&0\le t\le1/2,\\
2(1-t)^2,&1/2<t\le1,\\
0,&t>1.
\end{cases}
\end{equation}
Set $h_{\varnothing}:=0$.

\begin{lemma}[Convex Smoothing]
The function $h_C$ belongs to $C^{1,1}(\mathbb{R}^d)$, takes values in $[0,1]$, and has convex positive superlevel sets. With $\mathcal S_C:=(C+MB_2^d)\setminus\operatorname{int}C$, its weak derivatives satisfy
\begin{equation}
\label{eq.test_derivatives}
\begin{aligned}
|Dh_C(x)[w]|&\le2\|M^{-1}w\|_2\mathbf1_{\mathcal S_C}(x),\\
|D^2h_C(x)[v,w]|&\le4\|M^{-1}v\|_2\|M^{-1}w\|_2
\mathbf1_{\mathcal S_C}(x)
\end{aligned}
\end{equation}
almost everywhere. In particular, $\|M^{-1}w\|_2\le\delta_j^{-1}\|w\|_2$ for $w\in E_{j-1}^{\perp}$. Moreover,
\begin{equation}
\label{eq.smoothing_inequality}
\Delta_{n,d}
\le\sup_{C\text{ closed convex}}|\mathbb{E}h_C(S)-\mathbb{E}h_C(Z)|
+c_{\mathrm s}L\beta.
\end{equation}
\end{lemma}

\begin{proof}
Put $A:=M^{-1}C$ and $r(y):=\operatorname{dist}(y,A)$, and let $p(y)$ be the metric projection onto $A$. Outside $A$, set $v(y):=(y-p(y))/r(y)$. Then
$$
\nabla r(y)=v(y),\qquad \|v(y)\|_2=1.
$$
The residual map $y\mapsto y-p(y)=\nabla(r^2/2)(y)$ is nonexpansive. Since $\psi'(t)/t$ is locally Lipschitz for $t>0$, the identity $\nabla(\psi\circ r)=(\psi'(r)/r)(y-p(y))$ gives a locally Lipschitz gradient outside $A$. On $r\le1/2$, this gradient equals $-4(y-p(y))$, so it is locally Lipschitz everywhere. Convexity of $r$ and differentiation of $\|\nabla r\|_2^2=1$ imply, almost everywhere outside $A$,
$$
0\preceq D^2r(y)\preceq r(y)^{-1}I_d,
\qquad
D^2r(y)v(y)=0.
$$
Thus the two terms in
$$
D^2(\psi\circ r)(y)
=\psi''(r(y))v(y)v(y)^\top+\psi'(r(y))D^2r(y)
$$
act on orthogonal subspaces. Since $|\psi''|\le4$ almost everywhere and $|\psi'(t)|/t\le4$ for $t>0$,
$$
\|D^2(\psi\circ r)(y)\|_{\mathrm{op}}
\le\max\left\{|\psi''(r(y))|,\frac{|\psi'(r(y))|}{r(y)}\right\}
\le4.
$$
Also $|\psi'|\le2$. Mollifying the gradient and using the uniform weak Hessian bound shows that the gradient is globally $4$-Lipschitz. Composing with $M^{-1}$ proves \eqref{eq.test_derivatives}. The derivatives vanish off $\mathcal S_C$, and the positive superlevel sets are convex neighborhoods of $C$ in the metric induced by $M^{-1}$. Orthogonality gives
$$
\|M^{-1}w\|_2^2
=\sum_{\ell\ge j}\delta_\ell^{-2}\|P_\ell w\|_2^2
\le\delta_j^{-2}\|w\|_2^2
\qquad(w\in E_{j-1}^{\perp}).
$$

For closed convex $C$, the function $h_C$ is an upper approximation to $\mathbf1_C$, with Gaussian error bounded by \eqref{eq.smoothing_cost}. Put $C^-:=C\ominus MB_2^d$. If $C^-\ne\varnothing$, then $h_{C^-}\le\mathbf1_C$ and
$$
\gamma_d(C)-\mathbb{E}h_{C^-}(Z)
\le\gamma_d(C)-\gamma_d(C^-)
\le c_{\mathrm s}L\beta.
$$
If $C^-=\varnothing$, the same inner boundary estimate gives $\gamma_d(C)\le c_{\mathrm s}L\beta$, so the zero function supplies the lower approximation. For a convex Borel set, use its closure for the upper approximation. The lower approximation constructed from the closure is positive only in its interior. For a full-dimensional convex set, that interior is contained in the original set. A convex set with empty interior has Gaussian measure zero. These observations prove \eqref{eq.smoothing_inequality} for the full class $\mathcal C_d$.
\end{proof}

\subsection{The independent comparison block}

Fix $i\in I_a$ and define
\begin{equation}
\label{eq.free_block}
Y:=\sum_{\ell\in I_a\setminus\{i\}}X_\ell,
\qquad
B:=\operatorname{Cov}(Y)\succeq bI_d.
\end{equation}
The other two complete blocks remain independent of $Y$. Suppose that, for some $\mathcal D\ge0$,
\begin{equation}
\label{eq.free_comparison}
d_{\mathcal C}(Y,G)\le\mathcal D,\qquad G\sim N(0,B).
\end{equation}
Section~\ref{sec.main_proof} will provide this comparison by induction. All the following bounds are uniform in deterministic shifts and remain valid after conditioning on the other two blocks and their resampling variables.

Define the Gaussian interpolation
\begin{equation}
U_th(x):=\mathbb{E}h\left(\frac{x+tZ}{\sqrt{1+t^2}}\right),
\qquad t>0.
\end{equation}
For the standard Gaussian density $\varphi_d$, define
$$
H_m(z)[v_1,\ldots,v_m]
:=(-1)^m\varphi_d(z)^{-1}D^m\varphi_d(z)[v_1,\ldots,v_m].
$$
The generating-function identity
$$
\mathbb{E}\left[
e^{\langle Z,s\rangle-\|s\|_2^2/2}
e^{\langle Z,t\rangle-\|t\|_2^2/2}\right]
=e^{\langle s,t\rangle}
$$
implies, by differentiating $m$ times in each variable,
$$
\mathbb{E}\bigl[H_m(Z)[v_1,\ldots,v_m]^2\bigr]
=\sum_{\sigma\in\mathfrak S_m}\prod_{r=1}^m
\langle v_r,v_{\sigma(r)}\rangle
\le m!\prod_{r=1}^m\|v_r\|_2^2.
$$
Consequently,
\begin{equation}
\label{eq.hermite_bounds}
\mathbb{E}|H_m(Z)[v_1,\ldots,v_m]|
\le\sqrt{m!}\prod_{r=1}^m\|v_r\|_2,
\qquad
\mathbb{E}|H_1(Z)[v]|
=\sqrt{\frac2\pi}\|v\|_2.
\end{equation}

\begin{lemma}[Gaussian Comparison]
Assume \eqref{eq.free_comparison}, and let $h:=h_C$. For deterministic $x,v_1,\ldots,v_m$ and $m\ge1$,
\begin{equation}
\label{eq.derivative_comparison}
\left|\mathbb{E}D^mU_th(Y+x)[v_1,\ldots,v_m]\right|
\le\sqrt{m!}\left[(b+t^2)^{-m/2}+\mathcal D t^{-m}\right]
\prod_{r=1}^m\|v_r\|_2.
\end{equation}
Define $\mathcal H:=c_{\mathrm s}b^{-1/2}L\beta+2\mathcal D$. Uniformly in deterministic $x,z$,
\begin{equation}
\label{eq.shell_comparison}
\mathbb{P}\left(\frac{Y+x+tz}{\sqrt{1+t^2}}\in\mathcal S_C\right)
\le\mathcal H\sqrt{1+t^2}.
\end{equation}
\end{lemma}

\begin{proof}
Every positive superlevel set of $h$ is convex. The layer-cake representation and \eqref{eq.free_comparison} therefore give
$$
\left|\mathbb{E}h\left(\frac{Y+x+tz}{\sqrt{1+t^2}}\right)
-\mathbb{E}h\left(\frac{G+x+tz}{\sqrt{1+t^2}}\right)\right|
\le\mathcal D.
$$
Denote the difference inside the absolute value by $g(z)$. Transferring all derivatives to the Gaussian kernel and using \eqref{eq.hermite_bounds},
$$
\begin{aligned}
&\left|\mathbb{E}D^mU_th(Y+x)[v_1,\ldots,v_m]
-\mathbb{E}D^mU_th(G+x)[v_1,\ldots,v_m]\right|\\
&\quad=t^{-m}\left|\mathbb{E}[H_m(Z)[v_1,\ldots,v_m]g(Z)]\right|\\
&\quad\le\sqrt{m!}\,\mathcal D t^{-m}\prod_{r=1}^m\|v_r\|_2.
\end{aligned}
$$
Put $V_t:=(B+t^2I_d)/(1+t^2)$. For the reference term,
$$
\begin{aligned}
&\left|\mathbb{E}D^mU_th(G+x)[v_1,\ldots,v_m]\right|\\
&\quad=(1+t^2)^{-m/2}
\left|\mathbb{E}\left[
h\left(V_t^{1/2}Z+\frac{x}{\sqrt{1+t^2}}\right)
H_m(Z)[V_t^{-1/2}v_1,\ldots,V_t^{-1/2}v_m]\right]\right|\\
&\quad\le\sqrt{m!}(1+t^2)^{-m/2}
\prod_{r=1}^m\|V_t^{-1/2}v_r\|_2\\
&\quad\le\sqrt{m!}(b+t^2)^{-m/2}\prod_{r=1}^m\|v_r\|_2.
\end{aligned}
$$
This proves \eqref{eq.derivative_comparison}.

The shell event is a difference of two convex events, so its comparison error is at most $2\mathcal D$. The Gaussian reference has covariance $B/(1+t^2)\succeq bI_d/(1+t^2)$. Lemma~\ref{lem.gaussian_sections} gives
$$
\begin{aligned}
\mathbb{P}\left(\frac{Y+x+tz}{\sqrt{1+t^2}}\in\mathcal S_C\right)
&\le\frac{c_{\mathrm s}L\beta}{\sqrt b}\sqrt{1+t^2}+2\mathcal D\\
&\le\mathcal H\sqrt{1+t^2}.
\end{aligned}
$$
Gaussian measure assigns zero mass to the boundary of a full-dimensional convex set. If $C$ has empty interior, its Gaussian measure is zero, and the same estimate applies.
\end{proof}

We will also transfer only part of a derivative to the Gaussian kernel. For $s=1,2$ and $r\ge1$,
\begin{equation}
\label{eq.partial_transfer}
\begin{aligned}
&D^{r+s}U_th(x)[v_1,\ldots,v_r,w_1,\ldots,w_s]\\
&\quad=\frac{t^{-r}}{(1+t^2)^{s/2}}
\mathbb{E}\left[
H_r(Z)[v_1,\ldots,v_r]
D^sh\left(\frac{x+tZ}{\sqrt{1+t^2}}\right)[w_1,\ldots,w_s]\right].
\end{aligned}
\end{equation}
This follows by differentiating the Gaussian convolution in the sense of distributions. For $t>0$, the Gaussian perturbation avoids the exceptional Lebesgue null sets almost surely, also after averaging over independent random shifts by Fubini. We may therefore set the weak derivatives to zero on those sets so that \eqref{eq.test_derivatives} holds everywhere. Conditional on $Z$, the two resampling blocks, and their proposals, the uniform shell bound \eqref{eq.shell_comparison} then applies. The factor $\sqrt{1+t^2}$ in that bound is absorbed by the factor $(1+t^2)^{-s/2}$ in \eqref{eq.partial_transfer}.

\section{Conditional projections and bounded resampling}
\label{sec.resampling}

In this section, we derive identities for expected derivatives by replacing summands with independent copies and accepting only replacements with bounded displacement. A conditional expectation estimate controls the random coefficients in these identities.

\begin{lemma}[Conditional Projections]
\label{lem.conditional_projections}
Let $V,V'$ be independent and identically distributed, and let $k(V,V')$ be a square-integrable random vector in a real Hilbert space. Suppose that $\mathbb{E}k(V,V')=0$ and $k(y,x)=\sigma k(x,y)$ for $\sigma\in\{-1,1\}$. If $m(V):=\mathbb{E}[k(V,V')\mid V]$, then
$$
\mathbb{E}\|m(V)\|^2\le\frac12\mathbb{E}\|k(V,V')\|^2.
$$
\end{lemma}

\begin{proof}
Symmetry gives $\mathbb{E}[k(V,V')\mid V']=\sigma m(V')$. Define
$$
r(V,V'):=k(V,V')-m(V)-\sigma m(V').
$$
Then $\mathbb{E}[r\mid V]=\mathbb{E}[r\mid V']=0$. Also $m(V)$ and $m(V')$ are independent and centered. All cross terms therefore vanish, and
$$
\begin{aligned}
\mathbb{E}\|k(V,V')\|^2
&=\mathbb{E}\|m(V)+\sigma m(V')\|^2+\mathbb{E}\|r(V,V')\|^2\\
&=2\mathbb{E}\|m(V)\|^2+\mathbb{E}\|r(V,V')\|^2
\ge2\mathbb{E}\|m(V)\|^2.
\end{aligned}
$$
\end{proof}

Fix a complete block $I_a$, a radius $R:=R_j$, and a direction $u\in E_j$. Suppress the block and radius subscripts by writing $N:=|I_a|$, $\Gamma:=\Gamma_a(R)$, and $c_u:=\Gamma^\dagger u$. Define
$$
\ell_i:=\mathbb{E}\bigl[-D_i\mathbf1_{\{\|D_i\|_2\le R\}}\mid X_i\bigr],
\qquad
A_i:=\frac12\mathbb{E}\bigl[D_iD_i^\top\mathbf1_{\{\|D_i\|_2\le R\}}\mid X_i\bigr],
$$
and set
$$
\ell:=\sum_{i\in I_a}\ell_i,\qquad
V:=\sum_{i\in I_a}(A_i-\mathbb{E}A_i),\qquad
s(u):=c_u^\top\ell,\qquad t(u):=Vc_u.
$$

Choose $J$ uniformly from $I_a$, independently of all vectors. Replace $X_J$ by $\widehat X_J$ if $\|\widehat X_J-X_J\|_2\le R$, and keep $X_J$ otherwise. The displacement of the block sum is
$$
\zeta:=D_J\mathbf1_{\{\|D_J\|_2\le R\}}.
$$
The acceptance event is invariant under interchanging $X_J$ and $\widehat X_J$. Thus the original and updated block configurations form an exchangeable pair. Conditional on the original block,
\begin{equation}
\label{eq.exchange_regression}
\mathbb{E}[\zeta\mid(X_i)_{i\in I_a}]=-\frac{\ell}{N},
\qquad
\frac N2\mathbb{E}[\zeta\zeta^\top\mid(X_i)_{i\in I_a}]=\Gamma+V.
\end{equation}

\begin{lemma}[Resampling Moments]
The coefficients above satisfy
\begin{equation}
\label{eq.score_moments}
\mathbb{E}s(u)^2\le\lambda^{-1}\|u\|_2^2,
\qquad
\mathbb{E}\|t(u)\|_2^2\le\frac{R^2}{4\lambda}\|u\|_2^2.
\end{equation}
Moreover,
\begin{equation}
\label{eq.remainder_moment}
N\mathbb{E}\bigl[|c_u^\top\zeta|\|\zeta\|_2^3\bigr]
\le4\lambda^{-1/2}\sqrt{\beta R^3}\,\|u\|_2.
\end{equation}
\end{lemma}

\begin{proof}
Put $e(u):=u^\top\Gamma^\dagger u=c_u^\top\Gamma c_u$. By \eqref{eq.inverse_energy}, $e(u)\le\lambda^{-1}\|u\|_2^2$. The scalar kernel $-c_u^\top D_i\mathbf1_{\{\|D_i\|_2\le R\}}$ is antisymmetric. Lemma~\ref{lem.conditional_projections} and independence across indices give
$$
\begin{aligned}
\mathbb{E}s(u)^2
&=\sum_{i\in I_a}\mathbb{E}(c_u^\top\ell_i)^2\\
&\le\frac12\sum_{i\in I_a}
\mathbb{E}\bigl[(c_u^\top D_i)^2\mathbf1_{\{\|D_i\|_2\le R\}}\bigr]\\
&=c_u^\top\Gamma c_u=e(u)
\le\lambda^{-1}\|u\|_2^2.
\end{aligned}
$$
For the covariance error, apply the symmetric case of the same lemma to the centered version of
$$
k_i(X_i,\widehat X_i):=\frac12D_iD_i^\top c_u\,
\mathbf1_{\{\|D_i\|_2\le R\}}.
$$
Its conditional mean is $A_ic_u$. Consequently,
$$
\begin{aligned}
\mathbb{E}\|t(u)\|_2^2
&=\sum_{i\in I_a}\mathbb{E}\|(A_i-\mathbb{E}A_i)c_u\|_2^2\\
&\le\frac12\sum_{i\in I_a}
\mathbb{E}\|k_i-\mathbb{E}k_i\|_2^2\\
&\le\frac18\sum_{i\in I_a}
\mathbb{E}\bigl[\|D_i\|_2^2(c_u^\top D_i)^2
\mathbf1_{\{\|D_i\|_2\le R\}}\bigr]\\
&\le\frac{R^2}{8}\sum_{i\in I_a}
\mathbb{E}\bigl[(c_u^\top D_i)^2\mathbf1_{\{\|D_i\|_2\le R\}}\bigr]\\
&=\frac{R^2}{4}e(u)
\le\frac{R^2}{4\lambda}\|u\|_2^2.
\end{aligned}
$$
Finally,
$$
N\mathbb{E}(c_u^\top\zeta)^2=2e(u),
$$
and \eqref{eq.displacement_moment} gives
$$
\begin{aligned}
N\mathbb{E}\|\zeta\|_2^4
&=\sum_{i\in I_a}\mathbb{E}\bigl[
\|D_i\|_2^4\mathbf1_{\{\|D_i\|_2\le R\}}\bigr]\\
&\le R\sum_{i\in I_a}\mathbb{E}\|D_i\|_2^3
\le8R\beta.
\end{aligned}
$$
As $\|\zeta\|_2\le R$, Cauchy--Schwarz yields
$$
\begin{aligned}
N\mathbb{E}\bigl[|c_u^\top\zeta|\|\zeta\|_2^3\bigr]
&\le R\sqrt{N\mathbb{E}(c_u^\top\zeta)^2}
\sqrt{N\mathbb{E}\|\zeta\|_2^4}\\
&\le4\sqrt{e(u)}\,\sqrt{\beta R^3}
\le4\lambda^{-1/2}\sqrt{\beta R^3}\,\|u\|_2.
\end{aligned}
$$
\end{proof}

\begin{lemma}[Resampling Identity]
Let $W$ be the sum of the block $I_a$, an independent random vector, and a deterministic shift. For a bounded $C^3$ function $g$ with bounded derivatives,
\begin{equation}
\label{eq.exchange_identity}
\begin{aligned}
\mathbb{E}D_ug(W)
={}&\mathbb{E}[s(u)g(W)]-\mathbb{E}D_{t(u)}g(W)\\
&+\frac N4\mathbb{E}\left[(c_u^\top\zeta)
\int_0^1r(1-r)D^3g(W+r\zeta)[\zeta,\zeta,\zeta]\,dr\right].
\end{aligned}
\end{equation}
\end{lemma}

\begin{proof}
Write $W':=W+\zeta$. Exchangeability and \eqref{eq.exchange_regression} imply
$$
\begin{aligned}
\mathbb{E}[(c_u^\top\zeta)(g(W')-g(W))]
&=-2\mathbb{E}[(c_u^\top\zeta)g(W)]\\
&=\frac2N\mathbb{E}[s(u)g(W)].
\end{aligned}
$$
For a $C^3$ function $F$ on $[0,1]$, two integrations by parts give
$$
\begin{aligned}
\int_0^1r(1-r)F'''(r)\,dr
&=-\int_0^1(1-2r)F''(r)\,dr\\
&=F'(1)+F'(0)-2(F(1)-F(0)).
\end{aligned}
$$
Apply this identity to $F(r):=g(W+r\zeta)$. After multiplication by $c_u^\top\zeta$, exchangeability identifies the two endpoint terms. Therefore
$$
\begin{aligned}
\mathbb{E}[s(u)g(W)]
={}&\frac N2\mathbb{E}[(c_u^\top\zeta)D_\zeta g(W)]\\
&-\frac N4\mathbb{E}\left[(c_u^\top\zeta)
\int_0^1r(1-r)D^3g(W+r\zeta)[\zeta,\zeta,\zeta]\,dr\right]\\
={}&\mathbb{E}D_{\Gamma c_u+Vc_u}g(W)\\
&-\frac N4\mathbb{E}\left[(c_u^\top\zeta)
\int_0^1r(1-r)D^3g(W+r\zeta)[\zeta,\zeta,\zeta]\,dr\right].
\end{aligned}
$$
Since $\Gamma c_u=u$ and $Vc_u=t(u)$, rearrangement gives \eqref{eq.exchange_identity}.
\end{proof}

\subsection{Two independent resamplings}

In the setting of \eqref{eq.free_block}, label the other two complete blocks by $1$ and $2$. Let $W$ be their sum, the free block $Y$, and a deterministic shift. Use the same radius $R_j$ for both resamplings, and denote their coefficients by $s_\ell,t_\ell,\zeta_\ell,N_\ell,c_{\ell,u}$ for $\ell=1,2$. For $u,v\in E_j$ and any $w$, apply \eqref{eq.exchange_identity} first in direction $u$, and then to the remaining main term in direction $v$. For a bounded $C^5$ function $g$ with bounded derivatives, this gives
\begin{equation}
\label{eq.two_block_identity}
\begin{aligned}
\mathbb{E}D^3g(W)[u,v,w]
={}&\mathbb{E}[s_1(u)s_2(v)D_wg(W)]\\
&-\mathbb{E}D^3g(W)[t_1(u),v,w]\\
&-\mathbb{E}[s_1(u)D^2g(W)[t_2(v),w]]\\
&+\frac{N_1}{4}\mathbb{E}\left[(c_{1,u}^\top\zeta_1)
\int_0^1r(1-r)D^5g(W+r\zeta_1)
[v,w,\zeta_1,\zeta_1,\zeta_1]\,dr\right]\\
&+\frac{N_2}{4}\mathbb{E}\left[s_1(u)(c_{2,v}^\top\zeta_2)
\int_0^1r(1-r)D^4g(W+r\zeta_2)
[w,\zeta_2,\zeta_2,\zeta_2]\,dr\right].
\end{aligned}
\end{equation}
In the second application, condition on the first block so that $s_1(u)$ is fixed. In every term of \eqref{eq.two_block_identity}, conditioning on the two complete blocks and their proposals leaves the law of $Y$ unchanged.

Set $\rho_j:=\sqrt{\beta R_j^3}$. After factoring out the deterministic direction norms, the five terms in \eqref{eq.two_block_identity} have derivative orders and absolute weights
\begin{equation}
\label{eq.five_weights}
\begin{array}{c|ccccc}
\text{derivative order}&1&3&2&5&4\\ \hline
\text{weight}&c&cR_j&cR_j&c\rho_j&c\rho_j.
\end{array}
\end{equation}
Indeed, the first three bounds follow from \eqref{eq.score_moments} and independence:
$$
\begin{aligned}
\mathbb{E}|s_1(u)s_2(v)|&\le\lambda^{-1}\|u\|_2\|v\|_2,\\
\mathbb{E}\|t_1(u)\|_2&\le\frac{R_j}{2\sqrt\lambda}\|u\|_2,\\
\mathbb{E}[|s_1(u)|\|t_2(v)\|_2]&\le\frac{R_j}{2\lambda}\|u\|_2\|v\|_2.
\end{aligned}
$$
For the remainders, integrate the kernel $r(1-r)$ and apply \eqref{eq.remainder_moment}. The fourth weight is bounded by
$$
\frac14\int_0^1r(1-r)\,dr\cdot4\lambda^{-1/2}\rho_j
=\frac{\rho_j}{6\sqrt\lambda}.
$$
For the fifth weight, independence supplies the additional factor $\mathbb{E}|s_1(u)|\le\lambda^{-1/2}\|u\|_2$. All five coefficients in \eqref{eq.five_weights} are therefore universal.

\section{Cubic energy and the Stein integral}
\label{sec.stein_integral}

We now combine the estimates from Sections~\ref{sec.smoothing} and~\ref{sec.resampling} to prove the integrated derivative bound in Proposition~\ref{prop.integrated_derivatives}. Each layer term must have its first two directions in $E_j$ for the two resamplings and its last two in $E_{j-1}^{\perp}$ for the Hessian bound at width $\delta_j$. The cubic decomposition below meets both requirements with the middle direction in $H_j$, and its energy estimate allows us to sum the layer contributions by Cauchy--Schwarz.

Continue to use $W$ consisting of two complete blocks, the free block $Y$, and a deterministic shift. Assume \eqref{eq.free_comparison}, $\beta\le\beta_0$, and $L\beta\le1$. For $h:=h_C$, define the symmetric trilinear form
$$
T_t(u,v,w):=\mathbb{E}D^3U_th(W)[u,v,w].
$$

\subsection{An orthogonal cubic decomposition}

\begin{lemma}[Cubic Energy]
Let $T$ be a symmetric trilinear form and $u\in\mathbb{R}^d$. Put $p_j:=P_ju$ and $x_j(r):=\Pi_{j-1}u+rp_j$. For every integer $k\le0$,
\begin{equation}
\label{eq.cubic_decomposition}
T(u,u,u)
=T(\Pi_ku,\Pi_ku,3u-2\Pi_ku)
+6\sum_{j=k+1}^0\int_0^1
T(x_j(r),p_j,u-x_j(r))\,dr.
\end{equation}
The early term satisfies
\begin{equation}
\label{eq.early_direction_norm}
\|\Pi_ku\|_2^2\|3u-2\Pi_ku\|_2
\le3\|u\|_2^3.
\end{equation}
For every nonnegative square-summable sequence $(c_j)_{j\le0}$,
\begin{equation}
\label{eq.cubic_energy_bound}
6\sum_{j\le0}c_j\int_0^1
\|x_j(r)\|_2\|p_j\|_2\|u-x_j(r)\|_2\,dr
\le\sqrt6\,\|u\|_2^3\left(\sum_{j\le0}c_j^2\right)^{1/2}.
\end{equation}
In the $j$th integrand of \eqref{eq.cubic_decomposition}, the directions lie in $E_j$, $H_j$, and $E_{j-1}^{\perp}$, respectively.
\end{lemma}

\begin{proof}
Define $F(x):=3T(x,x,u)-2T(x,x,x)$. Along $x=x_j(r)$,
$$
\begin{aligned}
\frac{d}{dr}F(x_j(r))
&=6T(x_j(r),p_j,u)-6T(x_j(r),p_j,x_j(r))\\
&=6T(x_j(r),p_j,u-x_j(r)).
\end{aligned}
$$
As $\Pi_0u=u$, telescoping gives
$$
\begin{aligned}
T(u,u,u)
&=F(u)
=F(\Pi_ku)+\sum_{j=k+1}^0\bigl(F(\Pi_ju)-F(\Pi_{j-1}u)\bigr)\\
&=T(\Pi_ku,\Pi_ku,3u-2\Pi_ku)
+6\sum_{j=k+1}^0\int_0^1T(x_j(r),p_j,u-x_j(r))\,dr.
\end{aligned}
$$
The subspace assertions follow directly from the orthogonal decomposition.

By homogeneity, it suffices to prove the norm estimates when $\|u\|_2=1$. If $a:=\|\Pi_ku\|_2^2$, orthogonality gives
$$
\|\Pi_ku\|_2^2\|3u-2\Pi_ku\|_2
=a\sqrt{a+9(1-a)}
=a\sqrt{9-8a}
\le3.
$$

To prove \eqref{eq.cubic_energy_bound}, put $a_j:=\|p_j\|_2^2$, so $\sum_j a_j=1$. Orthogonality and direct integration yield
$$
\begin{aligned}
&\int_0^1\|x_j(r)\|_2^2\|p_j\|_2^2\|u-x_j(r)\|_2^2\,dr\\
&\quad=a_j\int_0^1
\left(\sum_{\ell<j}a_\ell+r^2a_j\right)
\left(\sum_{\ell>j}a_\ell+(1-r)^2a_j\right)\,dr\\
&\quad=a_j\left(\sum_{\ell<j}a_\ell\right)
\left(\sum_{\ell>j}a_\ell\right)
+\frac{a_j^2(1-a_j)}3+\frac{a_j^3}{30}.
\end{aligned}
$$
The first terms sum to the third elementary symmetric polynomial:
$$
\sum_j a_j\left(\sum_{\ell<j}a_\ell\right)
\left(\sum_{\ell>j}a_\ell\right)
=\sum_{i<j<k}a_ia_ja_k
=\frac{1-3\sum_j a_j^2+2\sum_j a_j^3}{6}.
$$
Consequently,
\begin{equation}
\begin{aligned}
\sum_j\int_0^1\|x_j(r)\|_2^2\|p_j\|_2^2\|u-x_j(r)\|_2^2\,dr
&=\frac16-\frac16\sum_j a_j^2+\frac1{30}\sum_j a_j^3\\
&\le\frac16.
\end{aligned}
\end{equation}
The last inequality uses $0\le a_j\le1$. Cauchy--Schwarz on the product of counting measure and Lebesgue measure now gives
$$
\begin{aligned}
&6\sum_jc_j\int_0^1
\|x_j(r)\|_2\|p_j\|_2\|u-x_j(r)\|_2\,dr\\
&\quad\le6\left(\sum_jc_j^2\right)^{1/2}
\left(\sum_j\int_0^1
\|x_j(r)\|_2^2\|p_j\|_2^2\|u-x_j(r)\|_2^2\,dr\right)^{1/2}\\
&\quad\le\sqrt6\left(\sum_jc_j^2\right)^{1/2}.
\end{aligned}
$$
Rescaling $u$ proves \eqref{eq.cubic_energy_bound}.
\end{proof}

For each $t>0$, choose
\begin{equation}
k(t):=\max\{j\le0:\delta_j\le t\}.
\end{equation}
This maximum exists because $\delta_j\to0$ as $j\to-\infty$. Here early and late refer to the layers with $\delta_j\le t$ and $t<\delta_j$, respectively. If $E_{k(t)}=\{0\}$, the early term in \eqref{eq.cubic_decomposition} vanishes. The identity is applied separately at each $t$, without differentiating $k(t)$ or its projection.

\subsection{The late layers}

\begin{lemma}[Late Layers]
\label{lem.late_layers}
For unit vectors $u\in E_j$, $v\in H_j$, and $w\in E_{j-1}^{\perp}$,
\begin{equation}
\label{eq.late_layer_integral}
\int_0^{\delta_j}|T_t(u,v,w)|\frac{t}{1+t^2}\,dt
\le c\left[\delta_j+\frac{R_j}{\delta_j}+R_j
+\frac{\mathcal H\sqrt{\rho_j}}{\delta_j^2}\right].
\end{equation}
The total integral of the absolute values of the late terms in \eqref{eq.cubic_decomposition} is at most
$$
c\left[L\beta+\frac1L+\beta+\frac{\mathcal H}{L^2\beta}\right]\|u\|_2^3.
$$
\end{lemma}

\begin{proof}
For the first assertion, write $R:=R_j$, $\delta:=\delta_j$, and $\rho:=\rho_j$. Retain the Hessian in directions $v,w$ and transfer the $u$ derivative to the Gaussian kernel. Equations \eqref{eq.test_derivatives}, \eqref{eq.hermite_bounds}, \eqref{eq.shell_comparison}, and \eqref{eq.partial_transfer} imply
\begin{equation}
\label{eq.direct_late_bound}
\begin{aligned}
\frac{t}{1+t^2}|T_t(u,v,w)|
&\le\frac{4}{\delta^2(1+t^2)^2}
\sqrt{\frac2\pi}\,\mathcal H\sqrt{1+t^2}\\
&\le\frac{c\mathcal H}{\delta^2}.
\end{aligned}
\end{equation}

For a second estimate, apply \eqref{eq.two_block_identity} to $g:=U_th$. Its main term has one derivative, so its weighted absolute value is at most $c/(1+t^2)$. For the first covariance term, either retain the Hessian in directions $v,w$ or transfer all derivatives to the Gaussian kernel. For the second covariance term, either retain the gradient in direction $w$ or transfer both derivatives. After multiplying by $t/(1+t^2)$, these estimates give, respectively,
$$
cR\min\{\delta^{-2},t^{-2}\},
\qquad
cR\min\{\delta^{-1},t^{-1}\}.
$$
Here the global derivative bounds suffice, so no shell factor is used.

In the fifth-order remainder, retain the Hessian in directions $v,w$ and transfer the three displacement derivatives. Conditional on both resampling blocks and their proposals, the shell bound applies uniformly in the interpolation parameter. Its weighted contribution is at most
$$
\frac{c\rho}{t^2(1+t^2)^2\delta^2}\,
\mathcal H\sqrt{1+t^2}
\le\frac{c\rho\mathcal H}{t^2\delta^2}.
$$
For the fourth-order remainder, retain the gradient in direction $w$ and transfer the three displacement derivatives. Its weighted contribution is at most
$$
\frac{c\rho}{t^2(1+t^2)^{3/2}\delta}\,
\mathcal H\sqrt{1+t^2}
\le\frac{c\rho\mathcal H}{t^2\delta^2},
$$
where $\delta\le L\beta\le1$. Using $\min\{A,B+C\}\le B+\min\{A,C\}$ for $A,B,C\ge0$, we retain the main and covariance bounds and combine the remainder bound with \eqref{eq.direct_late_bound} to obtain
$$
\frac{t}{1+t^2}|T_t(u,v,w)|
\le c\left[
1+R\min\{\delta^{-2},t^{-2}\}
+R\min\{\delta^{-1},t^{-1}\}
+\frac{\mathcal H}{\delta^2}\min\{1,\rho t^{-2}\}\right].
$$
Integrating over $0<t<\delta$ gives
$$
\begin{aligned}
\int_0^\delta|T_t(u,v,w)|\frac{t}{1+t^2}\,dt
&\le c\left[
\delta+\frac R\delta+R
+\frac{\mathcal H}{\delta^2}\int_0^\infty\min\{1,\rho t^{-2}\}\,dt\right]\\
&=c\left[\delta+\frac R\delta+R+\frac{2\mathcal H\sqrt\rho}{\delta^2}\right].
\end{aligned}
$$
This proves \eqref{eq.late_layer_integral}.

Use homogeneity for the directions in \eqref{eq.cubic_decomposition}, then apply \eqref{eq.cubic_energy_bound} to the integrated coefficients. The triangle inequality in $\ell^2$ reduces the result to four geometric sums:
\begin{equation}
\label{eq.layer_square_sums}
\begin{aligned}
\left(\sum_{j\le0}\delta_j^2\right)^{1/2}
&=\frac{L\beta}{\sqrt{1-q^{-2/3}}},\\
\left(\sum_{j\le0}(R_j/\delta_j)^2\right)^{1/2}
&=\frac{r_0}{L\sqrt{1-q^{-4/3}}},\\
\left(\sum_{j\le0}R_j^2\right)^{1/2}
&=\frac{r_0\beta}{\sqrt{1-q^{-2}}},\\
\left(\sum_{j\le0}\frac{\rho_j}{\delta_j^4}\right)^{1/2}
&=\frac{r_0^{3/4}}{L^2\beta\sqrt{1-q^{-1/6}}}.
\end{aligned}
\end{equation}
Substitution proves the second assertion.
\end{proof}

\subsection{The early directions}

\begin{lemma}[Early Directions]
\label{lem.early_directions}
The early term in \eqref{eq.cubic_decomposition}, with $T:=T_t$ and $k:=k(t)$, satisfies
$$
\int_0^\infty
|T_t(\Pi_{k(t)}u,\Pi_{k(t)}u,3u-2\Pi_{k(t)}u)|
\frac{t}{1+t^2}\,dt
\le c\left[1+\frac{\mathcal D}{L}
+\frac{\mathcal D}{L^3\beta}\right]\|u\|_2^3.
$$
\end{lemma}

\begin{proof}
The first two directions lie in $E_{k(t)}$, so apply \eqref{eq.two_block_identity} at radius $R:=R_{k(t)}$. First take three unit directions. The main term is bounded by $c/(1+t^2)$ after multiplication by the integration weight. For each remaining term, condition on both complete blocks and use \eqref{eq.derivative_comparison}. With $\rho:=\rho_{k(t)}$, the four bounds are
\begin{equation}
\label{eq.early_four_terms}
\begin{aligned}
&cR\,\frac{t}{1+t^2}
\left[(b+t^2)^{-3/2}+\mathcal D t^{-3}\right],\\
&cR\,\frac{t}{1+t^2}
\left[(b+t^2)^{-1}+\mathcal D t^{-2}\right],\\
&c\rho\,\frac{t}{1+t^2}
\left[(b+t^2)^{-5/2}+\mathcal D t^{-5}\right],\\
&c\rho\,\frac{t}{1+t^2}
\left[(b+t^2)^{-2}+\mathcal D t^{-4}\right].
\end{aligned}
\end{equation}
The reference terms have finite integrals. In fact,
\begin{equation}
\label{eq.reference_integrals}
\begin{aligned}
\int_0^\infty\frac{t}{(1+t^2)(b+t^2)}\,dt
&=\frac12\int_0^\infty\frac{dx}{(1+x)(b+x)}
=\frac{\log(1/b)}{2(1-b)},\\
\int_0^\infty\frac{t}{1+t^2}(b+t^2)^{-m/2}\,dt
&\le\int_0^\infty t(b+t^2)^{-m/2}\,dt
=\frac{b^{1-m/2}}{m-2}\qquad(m>2).
\end{aligned}
\end{equation}
Together with $R\le r_0\beta$ and $\rho\le r_0^{3/2}\beta^2$, these estimates bound the reference contributions by a universal constant.

For the comparison errors, put $h_0:=L\beta$. The scale choice gives
\begin{equation}
R_{k(t)}\le r_0\beta\min\{1,(t/h_0)^3\},
\qquad
\rho_{k(t)}\le r_0^{3/2}\beta^2\min\{1,(t/h_0)^{9/2}\}.
\end{equation}
The required scalar integrals are
$$
\begin{aligned}
\int_0^\infty\min\{1,(t/h_0)^3\}t^{-2}\,dt
&=h_0^{-3}\int_0^{h_0}t\,dt+\int_{h_0}^\infty t^{-2}\,dt
=\frac3{2h_0},\\
\int_0^\infty\min\{1,(t/h_0)^{9/2}\}t^{-4}\,dt
&=h_0^{-9/2}\int_0^{h_0}t^{1/2}\,dt+\int_{h_0}^\infty t^{-4}\,dt
=\frac1{h_0^3}.
\end{aligned}
$$
Using $1/(1+t^2)\le1$ and $t/(1+t^2)\le1/2$, we obtain
\begin{equation}
\label{eq.early_error_integrals}
\begin{aligned}
\int_0^\infty\frac{R_{k(t)}}{t^2(1+t^2)}\,dt&\le\frac{3r_0}{2L},&
\int_0^\infty\frac{R_{k(t)}}{t(1+t^2)}\,dt&\le\frac{3r_0}{4L},\\
\int_0^\infty\frac{\rho_{k(t)}}{t^4(1+t^2)}\,dt
&\le\frac{r_0^{3/2}}{L^3\beta},&
\int_0^\infty\frac{\rho_{k(t)}}{t^3(1+t^2)}\,dt
&\le\frac{r_0^{3/2}}{2L^3\beta}.
\end{aligned}
\end{equation}
These bounds control all four errors in \eqref{eq.early_four_terms}. Finally, \eqref{eq.early_direction_norm} restores the original directions with a factor at most $3\|u\|_2^3$.
\end{proof}

\subsection{The integrated derivative estimate}

\begin{proposition}[Integrated Derivatives]
\label{prop.integrated_derivatives}
Assume $\beta\le\beta_0$, $L\ge1$, $L\beta\le1$, and \eqref{eq.free_comparison}. Uniformly in the deterministic shift in $W$,
\begin{equation}
\int_0^\infty
\left|\mathbb{E}D^3U_th(W)[u,u,u]\right|
\frac{t}{1+t^2}\,dt
\le\mathcal A\|u\|_2^3,
\end{equation}
where, for a universal constant $c$, one may take
\begin{equation}
\label{eq.integrated_constant}
\mathcal A:=c\left(1+\frac{\mathcal D}{L^2\beta}\right).
\end{equation}
\end{proposition}

\begin{proof}
Apply \eqref{eq.cubic_decomposition} at $k=k(t)$ for each $t>0$. The late terms are active only on $0<t<\delta_j$, so Lemma~\ref{lem.late_layers} applies after interchanging their finite sum with integration. Lemma~\ref{lem.early_directions} controls the early term over $(0,\infty)$. Together, the two estimates bound the integral, divided by $\|u\|_2^3$, by
$$
c\left[
1+L\beta+\frac1L+\beta+\frac{\mathcal H}{L^2\beta}
+\frac{\mathcal D}{L}+\frac{\mathcal D}{L^3\beta}
\right].
$$
Since $L\ge1$, $L\beta\le1$, and $\mathcal H=c_{\mathrm s}b^{-1/2}L\beta+2\mathcal D$,
$$
\begin{aligned}
\frac{\mathcal H}{L^2\beta}
&=\frac{c_{\mathrm s}b^{-1/2}}L
+\frac{2\mathcal D}{L^2\beta},\\
\frac{\mathcal D}{L}
&=(L\beta)\frac{\mathcal D}{L^2\beta}
\le\frac{\mathcal D}{L^2\beta},\\
\frac{\mathcal D}{L^3\beta}
&\le\frac{\mathcal D}{L^2\beta}.
\end{aligned}
$$
Absorbing the bounded terms into a universal constant proves \eqref{eq.integrated_constant}.
\end{proof}

\section{Proof of the main theorem}
\label{sec.main_proof}

In this section, we combine Proposition~\ref{prop.integrated_derivatives} with the Stein equation and the smoothing inequality \eqref{eq.smoothing_inequality} to prove Theorem~\ref{thm.main_result} by induction on the number of summands.

Fix $h:=h_C$ and define
\begin{equation}
\label{eq.stein_solution}
f(x):=-\int_0^\infty\bigl(U_th(x)-\mathbb{E}h(Z)\bigr)
\frac{t}{1+t^2}\,dt.
\end{equation}
The bounded gradient and bounded weak Hessian of $h$ justify differentiation through order three. In particular, transferring one derivative to the Gaussian kernel gives
$$
|D^3U_th(x)[u,v,w]|
\le\frac{\|D^2h\|_\infty}{t(1+t^2)}
\|u\|_2\|v\|_2\|w\|_2.
$$
Multiplication by $t/(1+t^2)$ leaves the integrable factor $(1+t^2)^{-2}$. The first two derivatives are bounded by the gradient and Hessian bounds with factors $(1+t^2)^{-1/2}$ and $(1+t^2)^{-1}$, respectively. Hence $f\in C^3$ and its first three derivatives are bounded. The integral defining $f$ converges at zero by boundedness of $h$. At infinity, couple $(x+tZ)/\sqrt{1+t^2}$ with $Z$ and use the Lipschitz bound on $h$.

Set $\tau:=\frac12\log(1+t^2)$. Then $d\tau=t(1+t^2)^{-1}\,dt$, and $U_t$ is the Ornstein--Uhlenbeck semigroup at time $\tau$. Its generator is $\mathcal Lf:=\operatorname{tr}D^2f-x^\top\nabla f$. Integrating the semigroup derivative in \eqref{eq.stein_solution} gives
\begin{equation}
\operatorname{tr}D^2f(x)-x^\top\nabla f(x)
=-\int_0^\infty\frac{d}{d\tau}
U_{\sqrt{e^{2\tau}-1}}h(x)\,d\tau
=h(x)-\mathbb{E}h(Z).
\end{equation}
One may first integrate over a finite time interval and then pass to the limit using the derivative bounds above.

For each $i$, put $S_i:=S-X_i$, and let $X_i'$ be an independent copy of $X_i$. Since $\sum_i\Sigma_i=I_d$,
$$
\begin{aligned}
\mathbb{E}\operatorname{tr}D^2f(S)
&=\sum_i\mathbb{E}D^2f(S)[X_i',X_i']\\
&=\sum_i\mathbb{E}D^2f(S_i)[X_i',X_i']
+\sum_i\int_0^1
\mathbb{E}D^3f(S_i+rX_i)[X_i',X_i',X_i]\,dr.
\end{aligned}
$$
Independence and centering of $X_i$ give
$$
\begin{aligned}
\mathbb{E}[X_i^\top\nabla f(S)]
&=\int_0^1\mathbb{E}D^2f(S_i+rX_i)[X_i,X_i]\,dr\\
&=\mathbb{E}D^2f(S_i)[X_i,X_i]
+\int_0^1\int_0^r
\mathbb{E}D^3f(S_i+sX_i)[X_i,X_i,X_i]\,ds\,dr\\
&=\mathbb{E}D^2f(S_i)[X_i,X_i]
+\int_0^1(1-r)
\mathbb{E}D^3f(S_i+rX_i)[X_i,X_i,X_i]\,dr.
\end{aligned}
$$
The second-order terms cancel because $X_i$ and $X_i'$ have the same covariance and are independent of $S_i$. Thus
\begin{equation}
\label{eq.outer_expansion}
\begin{aligned}
\mathbb{E}h(S)-\mathbb{E}h(Z)
={}&\sum_i\int_0^1
\mathbb{E}D^3f(S_i+rX_i)[X_i',X_i',X_i]\,dr\\
&-\sum_i\int_0^1(1-r)
\mathbb{E}D^3f(S_i+rX_i)[X_i,X_i,X_i]\,dr.
\end{aligned}
\end{equation}
All terms are integrable. Indeed,
$$
\mathbb{E}\bigl[\|X_i'\|_2^2\|X_i\|_2\bigr]
=\mathbb{E}\|X_i\|_2^2\,\mathbb{E}\|X_i\|_2
\le b_i^{2/3}b_i^{1/3}=b_i.
$$

Banach's theorem for symmetric multilinear forms on a real Hilbert space \cite{Banach1938} states that a symmetric trilinear form satisfies
\begin{equation}
\label{eq.banach_identity}
\sup_{\|a\|_2,\|b\|_2,\|c\|_2\le1}|T(a,b,c)|
=\sup_{\|u\|_2\le1}|T(u,u,u)|.
\end{equation}
For fixed $i$, $r\in[0,1]$, and $x\in\mathbb{R}^d$, define the symmetric trilinear form
$$
T_{i,r,x}(a,b,c):=\mathbb{E}D^3f(S_i+rx)[a,b,c].
$$
If $i\in I_a$, then $S_i+rx$ consists of the two other complete blocks, the free block in \eqref{eq.free_block}, and the fixed shift $rx$. The matrix $M$ depends on the distributions of the original array and is the same for every outer index. Whenever \eqref{eq.free_comparison} holds with the same $\mathcal D$ for all these free blocks, Proposition~\ref{prop.integrated_derivatives} and \eqref{eq.stein_solution} give, uniformly in $x$,
$$
|T_{i,r,x}(u,u,u)|\le\mathcal A\|u\|_2^3.
$$
The integrable derivative bounds above justify interchanging expectation with the Stein integral. By \eqref{eq.banach_identity}, $|T_{i,r,x}(a,b,c)|\le\mathcal A\|a\|_2\|b\|_2\|c\|_2$. Since $S_i$ is independent of $X_i,X_i'$, conditioning on these vectors fixes the shift and directions. Taking expectations gives
$$
\begin{aligned}
\left|\mathbb{E}D^3f(S_i+rX_i)[X_i',X_i',X_i]\right|
&\le\mathcal A\,\mathbb{E}\bigl[\|X_i'\|_2^2\|X_i\|_2\bigr]
\le\mathcal A b_i,\\
\left|\mathbb{E}D^3f(S_i+rX_i)[X_i,X_i,X_i]\right|
&\le\mathcal A\,\mathbb{E}\|X_i\|_2^3
=\mathcal A b_i.
\end{aligned}
$$
Substitution into \eqref{eq.outer_expansion} yields
\begin{equation}
\label{eq.outer_bound}
\begin{aligned}
|\mathbb{E}h(S)-\mathbb{E}h(Z)|
&\le\sum_i\left(\int_0^1\mathcal A b_i\,dr
+\int_0^1(1-r)\mathcal A b_i\,dr\right)\\
&=\frac32\mathcal A\beta.
\end{aligned}
\end{equation}

\begin{proof}[Proof of Theorem~\ref{thm.main_result}]
We prove the existence of a finite universal constant by strong induction on $n$, simultaneously for all dimensions. Let $L\ge1$ and $C_0\ge\max\{\beta_0^{-1},L\}$ be universal constants to be chosen below. For $n=1$, Lyapunov's inequality gives $\beta\ge(\mathbb{E}\|X_1\|_2^2)^{3/2}=d^{3/2}\ge1$, so the assertion follows from $\Delta_{n,d}\le1$. For general $n$, the same trivial bound settles $\beta\ge1/C_0$.

Suppose that the assertion holds for fewer than $n$ summands and that $\beta<1/C_0$. Then $\beta\le\beta_0$ and $L\beta\le1$, so all the preceding constructions apply. For every $i\in I_a$, the free block $Y$ in \eqref{eq.free_block} contains fewer than $n$ summands and has covariance $B\succeq bI_d$. Its whitened summands are independent, centered, and have covariance sum $I_d$. The induction hypothesis covers all such arrays with finite third moments, so their third-moment sum need not be at most $\beta_0$. It therefore gives
$$
\begin{aligned}
d_{\mathcal C}(Y,N(0,B))
&=d_{\mathcal C}(B^{-1/2}Y,Z)\\
&\le C_0\sum_{\ell\in I_a\setminus\{i\}}
\mathbb{E}\|B^{-1/2}X_\ell\|_2^3\\
&\le C_0\|B^{-1/2}\|_{\mathrm{op}}^3\beta
\le b^{-3/2}C_0\beta.
\end{aligned}
$$
Thus \eqref{eq.free_comparison} holds uniformly in $i$ with $\mathcal D:=b^{-3/2}C_0\beta$. Proposition~\ref{prop.integrated_derivatives} now yields
$$
\mathcal A
\le c\left(1+\frac{b^{-3/2}C_0}{L^2}\right)
\le c\left(1+\frac{C_0}{L^2}\right).
$$
Here we used the fact that $b$ is fixed. Combining \eqref{eq.smoothing_inequality} and \eqref{eq.outer_bound}, there is a universal constant $C_1\ge1$, independent of $L$ and $C_0$, such that
\begin{equation}
\label{eq.qualitative_closure}
\frac{\Delta_{n,d}}{\beta}
\le c_{\mathrm s}L+\frac32\mathcal A
\le C_1(L+1)+\frac{C_1 C_0}{L^2}.
\end{equation}
Choose $L$ first so that $L^2\ge2C_1$, and then choose
$$
C_0\ge\max\{\beta_0^{-1},L,2C_1(L+1)\}.
$$
Both choices are finite and independent of the array. Substitution into \eqref{eq.qualitative_closure} gives
$$
\frac{\Delta_{n,d}}{\beta}
\le\frac{C_0}{2}+\frac{C_0}{2}=C_0.
$$
With these fixed choices, the induction closes for every $n$ and $d$. Hence $\sup_{d\ge1}K_d<\infty$.
\end{proof}

\clearpage
\appendix
\section{An explicit universal constant}
\label{app.explicit_constant}

This appendix gives a numerical choice of the universal constant in Theorem~\ref{thm.main_result}. We specify the parameters in the preceding construction, retain the coefficients in its estimates, and close the induction with $C_0=980$.

\subsection{Parameters and preliminary estimates}

Choose
\begin{equation}
\label{eq.numeric_parameters}
b:=\frac5{16},\qquad \lambda:=\frac18,\qquad
r_0:=20,\qquad q:=32,\qquad \beta_0:=\frac1{980}.
\end{equation}
The inequalities $980^2>98^3$ and $\sqrt{2/3}<49/60$ give
$$
\beta_0^{2/3}\le\frac1{98},\qquad
\sqrt{\frac23}\,\beta_0^{2/3}
<\frac{49}{60}\frac1{98}=\frac1{120}.
$$
Thus the first condition in \eqref{eq.small_budget_conditions} holds because
$$
\frac13-\lambda-\frac4{r_0}
=\frac13-\frac18-\frac15=\frac1{120}.
$$
For the second condition,
$$
\left(1+\sqrt{\frac23}\right)\beta_0^{2/3}
\le\frac1{98}+\frac1{120}<\frac1{48}=\frac13-b.
$$
The covariance partition consequently satisfies
$$
B_a\succeq\frac{13}{40}I_d,\qquad
B_a-\Sigma_i\succeq\frac5{16}I_d.
$$
The resampling radii and smoothing widths become
$$
R_j=20\beta\,32^j,\qquad
\delta_j=L\beta\,32^{j/3},\qquad
\rho_j=\sqrt{\beta R_j^3}.
$$
By \eqref{eq.dimension_budget}, $\sum_{j\le0}32^j\dim H_j\le1024/31$.

For the Gaussian boundary estimate, \cite[Theorem~1.2 and Proposition~1.1]{Raic2019} give, in dimension $r\ge1$,
$$
\begin{aligned}
\max\{\gamma_r(A+tB_2^r)-\gamma_r(A),
\gamma_r(A)-\gamma_r(A\ominus tB_2^r)\}
&\le t\left[\sqrt{\frac2\pi}+0.59(r^{1/4}-1)\right]\\
&\le\frac45\,tr^{1/4}.
\end{aligned}
$$
The section argument and \eqref{eq.smoothing_cost} therefore give a smoothing error bounded by
$$
\frac{\frac45(1024/31)^{1/4}}{(1-32^{-1/9})^{3/4}}L\beta
<4.52L\beta.
$$
We may take $c_{\mathrm s}:=4.52$. Since $4.52/\sqrt b<8.1$, the shell budget satisfies
\begin{equation}
\label{eq.numeric_shell}
\mathcal H=c_{\mathrm s}b^{-1/2}L\beta+2\mathcal D
\le8.1L\beta+2\mathcal D.
\end{equation}

Substituting $\lambda=1/8$ into the resampling moment estimates gives
$$
\mathbb{E}s(u)^2\le8\|u\|_2^2,\qquad
\mathbb{E}\|t(u)\|_2^2\le2R^2\|u\|_2^2,
$$
and
$$
N\mathbb{E}[|c_u^\top\zeta|\|\zeta\|_2^3]
\le8\sqrt2\,\sqrt{\beta R^3}\,\|u\|_2.
$$
Accordingly, the five weights in \eqref{eq.five_weights} may be taken to be
\begin{equation}
\label{eq.explicit_weights}
\begin{array}{c|ccccc}
\text{derivative order}&1&3&2&5&4\\ \hline
\text{weight}&8&\sqrt2R_j&4R_j&(\sqrt2/3)\rho_j&(4/3)\rho_j.
\end{array}
\end{equation}
The first three coefficients are $\lambda^{-1}$, $(2\sqrt\lambda)^{-1}$, and $(2\lambda)^{-1}$. The fourth is $(6\sqrt\lambda)^{-1}=\sqrt2/3$, and multiplication by $\lambda^{-1/2}$ gives the fifth coefficient $4/3$.

Throughout the next two subsections, assume $L\ge1$, $\beta\le1/980$, and $L\beta\le7/50$.

\subsection{Late-layer coefficients}

Fix a layer and write $R:=R_j$, $\delta:=\delta_j$, and $\rho:=\rho_j$. Take unit directions as in Lemma~\ref{lem.late_layers}. The direct estimate in \eqref{eq.direct_late_bound}, using $\mathbb{E}|H_1(Z)[u]|=\sqrt{2/\pi}<4/5$, gives
$$
\frac{t}{1+t^2}|T_t(u,v,w)|
\le\frac{4\sqrt{2/\pi}\,\mathcal H}{\delta^2(1+t^2)^{3/2}}
\le\frac{3.2\mathcal H}{\delta^2}.
$$

For the five terms in \eqref{eq.two_block_identity}, use \eqref{eq.explicit_weights}. Transferring the single derivative in the main term to the Gaussian kernel gives the weighted bound $8\sqrt{2/\pi}/(1+t^2)\le6.4/(1+t^2)$. The first covariance term is bounded by
$$
\sqrt2R\min\{3.2\delta^{-2},\sqrt6\,t^{-2}\}.
$$
Splitting the integral at $t=\delta\sqrt{\sqrt6/3.2}$ gives
$$
\begin{aligned}
\sqrt2R\int_0^\delta
\min\{3.2\delta^{-2},\sqrt6\,t^{-2}\}\,dt
&=\frac{\sqrt2R}{\delta}
\left(2\sqrt{3.2\sqrt6}-\sqrt6\right)\\
&<4.46\,\frac R\delta.
\end{aligned}
$$
The second covariance term is bounded by $4R\min\{1.6\delta^{-1},\sqrt2\,t^{-1}\}$. Splitting at $t=\delta\sqrt2/1.6$ yields
$$
\begin{aligned}
4R\int_0^\delta\min\{1.6\delta^{-1},\sqrt2\,t^{-1}\}\,dt
&=4\sqrt2R\left[1+\log\left(\frac{1.6}{\sqrt2}\right)\right]\\
&<6.36R.
\end{aligned}
$$

For the fifth-order remainder, retain the Hessian and transfer the three displacement derivatives. Its weighted bound is
$$
\frac{\sqrt2}{3}\rho\cdot4\delta^{-2}\cdot\sqrt6\,
\frac{\mathcal H\sqrt{1+t^2}}{t^2(1+t^2)^2}
\le\frac{4\sqrt{12}}3\frac{\rho\mathcal H}{t^2\delta^2}.
$$
For the fourth-order remainder, retain the gradient. The bound is
$$
\frac43\rho\cdot2\delta^{-1}\cdot\sqrt6\,
\frac{\mathcal H\sqrt{1+t^2}}{t^2(1+t^2)^{3/2}}
\le\frac{8\sqrt6}{3}\frac{\rho\mathcal H}{t^2\delta}.
$$
Taking the minimum of the direct estimate and the full resampling estimate preserves the nonnegative main and covariance bounds. The remaining contribution is at most
$$
\frac{\mathcal H}{\delta^2}
\min\left\{3.2,\sqrt6\left(\frac{4\sqrt2}{3}+\frac{8\delta}{3}\right)
\rho t^{-2}\right\}.
$$
For $a,c>0$, splitting at $t=\sqrt{c/a}$ gives
$$
\int_0^\infty\min\{a,ct^{-2}\}\,dt
=a\sqrt{\frac ca}+c\int_{\sqrt{c/a}}^\infty t^{-2}\,dt
=2\sqrt{ac}.
$$
As $\delta\le7/50$, the integral of the remainder contribution is bounded by
$$
\begin{aligned}
\frac{2\mathcal H\sqrt\rho}{\delta^2}
\sqrt{3.2\sqrt6\left(\frac{4\sqrt2}{3}+\frac{8\delta}{3}\right)}
&\le\frac{2\mathcal H\sqrt\rho}{\delta^2}
\sqrt{3.2\sqrt6\left(\frac{4\sqrt2}{3}+\frac{28}{75}\right)}\\
&<8.42\,\frac{\mathcal H\sqrt\rho}{\delta^2}.
\end{aligned}
$$
Thus the integrated coefficient for one layer is at most
\begin{equation}
6.4\delta_j+4.46\frac{R_j}{\delta_j}
+6.36R_j+8.42\frac{\mathcal H\sqrt{\rho_j}}{\delta_j^2}.
\end{equation}

The geometric factors in \eqref{eq.layer_square_sums} satisfy
$$
\begin{aligned}
(1-32^{-2/3})^{-1/2}&<1.054,&
(1-32^{-4/3})^{-1/2}&<1.005,\\
(1-32^{-2})^{-1/2}&<1.001,&
\frac{20^{3/4}}{\sqrt{1-32^{-1/6}}}&<14.28.
\end{aligned}
$$
Apply \eqref{eq.cubic_energy_bound}, the triangle inequality in $\ell^2$, and $\sqrt6<2.45$. The total late contribution, divided by $\|u\|_2^3$, is at most
$$
\begin{aligned}
&2.45\left[
6.4\cdot1.054L\beta
+4.46\cdot1.005\frac{20}{L}
+6.36\cdot1.001\cdot20\beta
+8.42\cdot14.28\frac{\mathcal H}{L^2\beta}\right]\\
&\quad=16.52672L\beta+\frac{219.6327}{L}
+311.95164\beta+\frac{294.58212\mathcal H}{L^2\beta}\\
&\quad\le17L\beta+\frac{220}{L}+312\beta
+\frac{295\mathcal H}{L^2\beta}.
\end{aligned}
$$

\subsection{Early-direction coefficients}

We first refine the direction factor in \eqref{eq.early_direction_norm}. For $\|u\|_2=1$, put $\alpha:=\|\Pi_ku\|_2^2$. Orthogonality gives
$$
\|\Pi_ku\|_2^2\|3u-2\Pi_ku\|_2
=\alpha\sqrt{9-8\alpha}.
$$
The derivative of its square is $18\alpha-24\alpha^2=6\alpha(3-4\alpha)$, so the maximum over $0\le\alpha\le1$ is attained at $\alpha=3/4$. Consequently,
$$
\|\Pi_ku\|_2^2\|3u-2\Pi_ku\|_2
\le\frac34\sqrt3<1.3.
$$
The main term therefore contributes at most
$$
1.3\cdot6.4\int_0^\infty\frac{dt}{1+t^2}
=1.3\cdot6.4\frac\pi2.
$$
The coefficients of the four terms in \eqref{eq.early_four_terms}, obtained by multiplying \eqref{eq.explicit_weights} by the corresponding Hermite bounds, are $\sqrt{12}R$, $4\sqrt2R$, $(\sqrt{240}/3)\rho$, and $(4\sqrt{24}/3)\rho$. For their reference contributions, use $R_{k(t)}\le20/980$ and $\rho_{k(t)}\le20^{3/2}/980^2$, followed by \eqref{eq.reference_integrals}. Their sum with the main term is at most
$$
\begin{aligned}
1.3\Bigg[
6.4\frac\pi2
&+\sqrt{12}\frac{20}{980}\frac1{\sqrt b}
+4\sqrt2\frac{20}{980}\frac{\log(1/b)}{2(1-b)}\\
&+\frac{\sqrt{240}}3\frac{20^{3/2}}{980^2}\frac{b^{-3/2}}3
+\frac{4\sqrt{24}}3\frac{20^{3/2}}{980^2}\frac1{2b}
\Bigg]
<13.363<14.
\end{aligned}
$$
For the two covariance comparison errors, \eqref{eq.early_error_integrals} gives the scalar bounds $30/L$ and $15/L$. Their combined coefficient is
$$
1.3\left(30\sqrt{12}+15\cdot4\sqrt2\right)
<246
$$
in front of $\mathcal D/L$. For the two remainder comparison errors, the scalar bounds in \eqref{eq.early_error_integrals} are $20^{3/2}/(L^3\beta)$ and $20^{3/2}/(2L^3\beta)$. The combined coefficient is
$$
1.3\cdot20^{3/2}\left(\frac{\sqrt{240}}3+\frac{2\sqrt{24}}3\right)
<981
$$
in front of $\mathcal D/(L^3\beta)$. The total early contribution is therefore at most
$$
\left[14+\frac{246\mathcal D}{L}
+\frac{981\mathcal D}{L^3\beta}\right]\|u\|_2^3.
$$
Together with the late estimate, this shows that the integrated derivative estimate in Proposition~\ref{prop.integrated_derivatives} holds with
\begin{equation}
\label{eq.explicit_integrated_constant}
\begin{aligned}
\mathcal A:={}&14+17L\beta+\frac{220}{L}+312\beta
+\frac{295\mathcal H}{L^2\beta}\\
&+\frac{246\mathcal D}{L}
+\frac{981\mathcal D}{L^3\beta}.
\end{aligned}
\end{equation}

\subsection{Numerical closure}

We now repeat the induction of Section~\ref{sec.main_proof} with $C_0:=980$ and $L:=135$. The case $\beta\ge1/C_0$ follows from $\Delta_{n,d}\le1$. If $\beta<1/C_0$, the parameters \eqref{eq.numeric_parameters} apply, and $L\beta<135/980<7/50$. Under the induction hypothesis, the free-block comparison satisfies
$$
d_{\mathcal C}(Y,N(0,B))
\le b^{-3/2}C_0\beta
=\left(\frac{16}{5}\right)^{3/2}C_0\beta
<5.73C_0\beta.
$$
Take $\mathcal D:=5.73C_0\beta$. Equation~\eqref{eq.numeric_shell} then gives $\mathcal H\le(8.1L+11.46C_0)\beta$. For $\beta\le1/C_0$, \eqref{eq.explicit_integrated_constant} yields
$$
\begin{aligned}
\mathcal A
&\le14+\frac{17L}{C_0}+\frac{220}{L}+\frac{312}{C_0}
+\frac{295(8.1L+11.46C_0)}{L^2}\\
&\qquad+\frac{246\cdot5.73}{L}
+\frac{981\cdot5.73C_0}{L^3}\\
&=14+\frac{17\cdot135}{980}+\frac{220}{135}+\frac{312}{980}
+\frac{295(8.1\cdot135+11.46\cdot980)}{135^2}\\
&\qquad+\frac{246\cdot5.73}{135}
+\frac{981\cdot5.73\cdot980}{135^3}\\
&<230.459.
\end{aligned}
$$
Consequently,
$$
\frac{\Delta_{n,d}}{\beta}
\le4.52L+\frac32\mathcal A
<610.2+345.6885
=955.8885<980.
$$
This closes the induction with $C_0=980$ and proves the explicit assertion of Theorem~\ref{thm.main_result}.

\clearpage
\bibliographystyle{alpha}
\bibliography{ref}

\end{document}